\documentclass[11pt]{amsart}
\usepackage{amsmath, amssymb, mathrsfs}
\usepackage{mathpazo}
\usepackage{xcolor}
\def\lbr{\left\{}
\def\rbr{\right\}}

\def\mso{{\mathscr O}}

\def\me{{\mathbb  E}}

\def\mr{{\mathbb  R}}

\def\mp{{\mathbb  P}}

\def\diff{{\partial }}
\def\E{{\me }}
\def\rom{{\rm }}

\def\lnfrac#1#2{\raise.7ex \hbox{\Small $#1$}
  \kern-.15em/\kern-.15em  \lower.2ex \hbox{\Small $#2$}}

\theoremstyle{plain}
\newtheorem{theorem}{Theorem}[section]

\theoremstyle{definition}

\newtheorem{remark}[theorem]{Remark}

\numberwithin{equation}{section}

\begin{document}

\title[Bismut's Formula]{Bismut's Gradient Formula for Vector Bundles}
\author{Elton P. Hsu}
\address{Department of Mathematics\\
        Northwestern University \\ Evanston, IL 60208 USA}
\email{ehsu\@@math.northwestern.edu}
\author{Zhenan Wang}
\address{Department of Mathematics\\
      Northwestern University\\ Evanston, IL 60208 USA}
\email{zn\_wang\@@math.northwestern.edu}
\thanks{The research of the first author was supported in part by a Simons Foundation Collaboration Grant for Mathematicians and a research grant administered by the Chinese University of Science and Technology.}

\begin{abstract}
We prove a general Bismut's formula for the gradient of a class of smooth Wiener functionals over vector bundles of a compact Riemannian manifold. This general formula can be used repeatedly for obtaining probabilistic representation of higher order covariant derivatives of solutions of the heat equation similar to the classical Bismut's representation for the covariant gradient of the heat kernel.
\end{abstract}
\maketitle
\centerline{(preliminary version)}

\tableofcontents

\section{Introduction} Let $M$ be a compact Riemannian manifold of dimension $n$ and $\mp_x$ the law of Brownian motion on $M$ starting from $x$. The now classical Bismut's formula (\cite{Bismut}) is a probabilistic representation of the gradient of the heat semigroup
$$P_t f(x) = \me_xf(X_t) = \int_M p(t,x,y) f(y)\, dy,$$
where $p(t,x,y)$ is the transition density function of the Brownian $X$ on $M$ (the heat kernel).  Several interesting objects are involved in this formula. On the orthonormal frame bundle $\mso(M)$ of $M$ the scalarized Ricci curvature tensor is realized as an $O(n)$-invariant, End$(\mr^n)$-valued function $\text{Ric}: u\mapsto \text{Ric}_u\in \text{End}(\mr^n)$. Let $\lbr U_t\rbr$ be a horizontal lift of the Brownian motion $X = \lbr X_t\rbr$ to the frame bundle $\mso(M)$. The first object of interest is the multiplicative Feynman-Kac functional $M = \lbr M_t\rbr$ defined by
$$dM_t +\frac12 \text{Ric}_{U_t}M_t = 0, \quad M_0 = I.$$
The second object of interest is the so-called stochastic anti-development $W = \lbr W_t\rbr$ of the Brownian motion $M$. It is the euclidean Brownian motion which drives the stochastic differential equation for the horizontal Brownian motion $U$, namely,
$$dU_t = \sum_{i=1}^n H_i(U_t)\circ dW_t^i,$$
where $\lbr H_i\rbr$ are the canonical horizontal vector fields on $\mso(M)$. The classical Bismut's formula is given by
$$T\nabla\me_xf(X_T) = \me_x\left[f(X_T)\int_0^T M_s\, dW_s\right].$$
Since $u(t,x) = \me_x f(X_t)$ is the solution of the heat equation on $M$ with the initial function $u(0, x) = f(x)$, the above identity gives a  probabilistic representation of the covariant gradient of the solution.  This representation is equivalent to the following probabilistic representation of the gradient of the heat kernel:
$$T\nabla_x\log p(t,x,y) = \me_{x,y; T}\left[\int_0^t M_s\, dW_s\right],$$
where $\me_{x,y; T}$ is the expectation with respect to the law $\mp_{x,y; T}$ of a Brownian bridge from $x$ to $y$ in time $T$.  It can be regarded as a form of integration by parts formula because of the absence of the gradient operator on the right side; see Hsu~\cite{Hsu} for a detailed exposition for this formula.

Ever since its appearance 30 years ago, this explicit representation has found applications in many areas. Besides its obvious significance in stochastic analysis on manifolds, it has also played important roles in certain problems in financial mathematics involving computations of option greeks. One naturally considers possible generalizations of this beautiful formula in several directions. Two possible directions come to mind immediately: higher order derivatives and functions in more general spaces. It is an important observation due to Norris~\cite{Norris} that these two possible generalizations can be combined into one framework, in which the function $f$ is replaced by a smooth map on a vector bundle and the process $X$ is lifted to the vector bundle by introducing, in addition to a compatible connection on the vector bundle, also a vertical motion driven by $X$ itself.  In this way, the formula for higher derivative can simply be obtained by from the general formula by a judicious choice of the bundle map and the vertical motion.

Bismut~\cite{Bismut} derived his formula by a perturbation method in path space, which was very much in line with the techniques available at the time when Malliavin calculus was studied intensively. This method was followed by and large by Norris~\cite{Norris} in his generalization of Bismut's formula in the setting we have just mentioned. The technical difficulties involved in this method is tremendous. In fact, up to now Norris' paper has not been fully digested by the probabilistic community. In view of the new ideas and technical tools introduced during the intervening decades (mainly gradient-heat semigroup commutation relations and It\^o calculus for diffusion processes on manifolds) and importance of these results, it is highly appropriate that we revisit these results with the goal of finding a proof more in line with the current state of art of stochastic analysis on manifolds.  In addition, in our recent investigation of functional inequalities in path and loop spaces of a Riemannian manifold, we feel the need to extend the framework adopted in Norris~\cite{Norris} by considering maps between two different vector bundles instead a self-map on a single vector bundle. Thus the purpose of this current work is to use It\^o calculus to prove a probabilistic representation of the gradient of the solution of a class of heat equations between two vector bundles. The basic idea of the current approach is explained in Hsu~\cite{Hsu1} for the simplest case of a trivial vector bundle.

{\sc Acknowledgment.} The origin of the current work can be traced back to the early collaboration of the first author with Professors Marc Arnaudon (Bordeaux), Zhongmin Qian (Oxford) and Anton Thalmaier (Luxembourg) for the purpose of finding an alternative approach to Norris~\cite{Norris}.  Although this collaboration did not result in a joint publication, the discussion with these individuals helped clarify
many confusions and misunderstandings concerning the problem.  Their early contribution to the current project is hereby gratefully acknowledged.
The authors would also like to express their gratitude to the School of Mathematics of the University of Science and Technology of China.The School provided financial support for the authors' visit in the summer of 2014, during which a portion of the present work was carried out.

\section{Notations and the statement of the problem}

We work on a Riemannian manifold $M$ of dimension $m$.  Let $\tau_1: E_1\rightarrow M$ and $\pi_2: E_2\rightarrow M$ be two Riemannian vector bundles over $M$ of dimensions $m_1$ and $m_2$, respectively.  Let $x = \lbr x_t,\, t\ge 0\rbr$ be a Riemannian Brownian motion on $M$ and $W = \lbr W_t, \, t\ge0\rbr$ its anti-development, which is a standard Brownian motion in $\mr^m$ adapted to the filtration of $x$. The horizontal lift of $x$, i.e. the horizontal Brownian motion $U = \lbr U_t,\, t\ge 0\rbr$, is defined by
$$\diff U_t=H_i(U_t)\circ\diff W_t^i,$$
where $\{H_i\}$ are canonical horizontal vector fields on $\mathscr{F}(M)$.   The (stochastic) parallel transport along the path $x[0, t]$ is  $\tau_t = U_0 U_t^{-1}$.

Following Norris~\cite{Norris}, we introduce a diffusion process on $E_1$ driven by $x$ (or equivalently, by $W$) in the following manner,
$$\diff y_t=v_0\diff t+v_{1,i}\circ\diff W^i_t,$$
where $v_0,v_1$ are vector fields on $E_1$.  This process is a lift of the Brownian motion on $M$, i.e., $\pi_1 y_t = x_t$.  The precise interpretation of the equation for $y$ is
$$\diff\tau^{-1}_ty_t=\tau^{-1_t}v_0\diff t+\tau^{-1}_tv_1\circ\diff W_t.$$

Consider a smooth map between the two vector bundles $f:E_1\rightarrow E_2$ which respects the fibres. Then $\tau_T^{-1}f(x_T, y_T)$ is an $E_{2,x}$-valued random variable and the expectation $\E_{x,y}\tau_T^{-1}f(x_T,y_T)$ is an element in $E_{2,x}$ which varies smoothly
as a function of the initial point $(x_0, y_0) = (x, y)$.  The focus of this paper is to give a probabilistic representation of the covariant gradient  $\nabla\E_{x,y}\tau_T^{-1}f(x_T,y_T)$ (with respect to the $x$ variable).

For a section of a vector bundle, it will be much more convenient for our exposition if we deal with its scalarization, which is a euclidean space valued function on the frame bundle of that vector bundle.  We introduce the frame bundle $\mathscr{F}(M)=\cup_{x\in M}U_x$, where $$U=(U_0,U_1,U_2)=\begin{cases}
 \displaystyle U_0: & \mathbb{R}^m\rightarrow T_xM \\
 U_1: & \mathbb{R}^{n_1}\rightarrow E_{1,x} \\
 U_2: & \mathbb{R}^{n_2}\rightarrow E_{2,x}
\end{cases}.$$
Each $U_i$ is an isometry between its respective domain and target spaces.
The scalarization of the bundle function $f:\, E_1\rightarrow E_2$ is defined by
$$F: \mathscr{F}(M)\times\mathbb{R}^{n_1}\rightarrow \mathbb{R}^{n_2}$$
$$F(U,Y):=U_2^{-1}f(x, U_1Y).$$

Let $Y = \lbr Y_t, \, t\ge0\rbr$ be the scalarization of the process $y$. Then it is the solution of the stochastic differential equation
$$\diff Y_t=V_0(U_t,Y_t)\diff t+\sum_iV_{1,i}(U_t,Y_t)\circ\diff W_t^i,$$
where $V_0(U_t,Y_t):=U_{1,t}^{-1}v_0(x_t,U_{1,t}Y_t)$ and $V_{1,i}(U_t,Y_t):=U_{1,t}^{-1}v_{1,i}(x_t,U_{1,t}Y_t)$ are the scalarization of $v_0(x_t, y_t)$ and $v_{1, i}(x_t, y_t)$, respectively.

For ease of notation, we combine the two processes $U$ and $Y$ into a single process $Z_t = (U_t, Y_t)$. It is clear that
$$\diff Z_t=(0,V_0)\diff t+\sum_i(H_i,V_{1,i})\circ\diff W_t^i.$$

To further simplify the notation, we will write $(H_i,0)$, $(0, V_0)$, and $(0, V_{1, i})$ simly as  $H_i$, $V_0$, and $V_{1,i}$, respectively. whenever there is not possibility of confusion. Then the generator of the diffusion process $Z_t$ is
$$L=\frac{1}{2}\sum_i(H_i+V_{1,i})^2+V_0.$$
And $F(z,t)=\E_z\{F(Z_t)\}$ is the solution to equation
$$\frac{\partial F}{\partial t}=LF.$$
The goal of this paper is to find a probabilistic representation of the gradient of the solution to this heat equation over the vector bundles.
Our main object of interest $\nabla\E_{x,y}\tau_T^{-1}f(x_T,y_T),$ becomes $\nabla^HF(z,T)$. Here $\nabla^H$ is the horizontal gradient defined by
$$\nabla^HF=\sum_j(H_jF)e_j,$$
where $\{e_j\}_{j=1}^{n_1}$ is the canonical orthonormal basis in $\mathbb{R}^{n_1}$.

\section{Commutation of the gradient and the heat semigroup}

The first step towards establishing the probabilistic representation of the gradient of a solution to the heat equation is to commute the gradient operator with the heat semigroup.  In this way we pass the gradient operator under the expectation on the bundle map. This will be followed by an integration by parts argument to remove the gradient operator from the bundle map.

Consider the equation for $\nabla^HF$
\begin{align*}
\frac{\partial\nabla^HF}{\partial t} &= \nabla^HLF\\
                                     &= L\nabla^HF+[\nabla^H,L]F
\end{align*}
We have
$$L=\frac{1}{2}\sum_iH_i^2+\frac{1}{2}\sum_iV_{1,i}^2+\frac{1}{2}\sum_i\left(H_iV_{1,i}+V_{1,i}H_i\right)+V_0.$$
By further expand the commutator bracket, we obtain
$$\frac{\partial\nabla^HF}{\partial t} = L\nabla^HF-\frac{1}{2}\text{Ric}\nabla^HF+\text{\rom{(1)}}+\text{\rom{(2)}}+\text{\rom{(3)}}+\text{\rom{(4)}}+\text{\rom{(5)}}$$
where
\begin{align*}\tag{I}\label{I}
\text{\rom{(1)}}&=\sum_j\left(\sum_i-\Omega_{ji}^{(2)}(H_iF)+\sum_i(DH_iF)(\Omega_{ji}^{(1)}Y)\right) e_j\\
 &+\frac{1}{2}\sum_j\left(\sum_i-(H_i\Omega_{ji}^{(2)})F+\sum_i(DF)((H_i\Omega_{ji}^{(1)})Y)\right) e_j
\end{align*}
\begin{align*}\tag{II}\label{II}
\text{\rom{(2)}}&=\sum_j\left(\sum_i(D^2F)(H_jV_{1,i})(V_{1,i})\right)e_j\\
  &+\frac{1}{2}\sum_jDF\left(\sum_i(DH_jV_{1,i})(V_{1,i})+(DV_{1,i})(H_jV_{1,i})\right)e_j
\end{align*}
\begin{equation*}\tag{III}\label{III}\text{\rom{(3)}}=\sum_j\left(\sum_i(DH_iF)(H_jV_{1,i})\right)e_j+\frac{1}{2}\sum_jDF\left(\sum_iH_iH_jV_{1,i}\right)e_j\end{equation*}

\begin{align*}\tag{IV}\label{IV}
\text{\rom{(4)}}&=\sum_j\left(\sum_i(D^2F)(V_{1,i})(\Omega_{ji}^{(1)}Y)\right)e_j\\
&+\sum_j\left(\sum_i-\Omega_{ji}^{(2)}(DF)(V_{1,i})\right)e_j\\
&+\frac{1}{2}\sum_jDF\left(\sum_i(DV_{1,i})(\Omega_{ji}^{(1)}Y)\right)e_j
\end{align*}
\begin{equation*}\tag{V}\label{V}\text{\rom{(5)}}=\sum_j(DF)(H_jV_0)e_j\end{equation*}

In the above expansion, $-\frac{1}{2}\text{Ric}\nabla^HF+\text{\rom{(1)}}$ comes from $[\nabla^H,\sum H_i^2]$, $\text{\rom{(2)}}$ comes from $[\nabla^H,\sum V_i^2]$, then $\text{\rom{(3)}}+\text{\rom{(4)}}$ comes from $[\nabla^H,\sum H_iV_{1,i}+V_{1,i}H_i]$ and $\text{\rom{(5)}}$ comes from $[\nabla^H,V_0]$. We also group $\text{\rom{(3)}}+\text{\rom{(4)}}$ so that all terms in their( corresponding commutator with curvature go to $\text{\rom{(4)}}$ and the remaining terms go to $\text{\rom{(3)}}$. The reason for doing so is to emphasize the effect of the curvature of the bundle $E_1$ and $E_2$.

For clarification, we say a few words about the notations used here concerning the curvature tensors (see Kobayashi and Nomizu~\cite{KobayashiN}). Let $H_i$ be the canonical horizontal vector fields of the product bundle $TM\oplus E_1\oplus E_2$. Then the commutator $\Omega_{ij} = [H_j,H_i]$ is a vertical vector field on the same product bundle with three components, i.e.,
$$\Omega_{ij} = (\Omega_{ji}^{(0)},\Omega_{ji}^{(1)},\Omega_{ji}^{(2)})\in \mathfrak{o}(m)\times\mathfrak{o}(n_1)\times\mathfrak{o}(n_2),$$
where $\Omega_{ji}^{(l)}$ are the (scalarized) curvature tensors for $l=0$ (on the manifold $M$), 1 (for the vector bundle $E_1$), and 2 (for the vector bundle $E_2)$.

In the above computation, we have also adopted the following conventions. If $Q$ is a vector in $\mathbb{R}^{n_1}$, the notation $DQ$ denotes the row vector $(D_iQ)_{i=1}^{n_1}$; if $Q$ has also other coordinates, say it is a vector in $\mathbb{R}^{m}\times\mathbb{R}^{n_1}$, then $DQ$ will just be $0$ in $\mathbb{R}^{m}$ and $(D_iQ)_{i=1}^{n_1}$ in $\mathbb{R}^{n_1}$. We use $(D^2Q)(\alpha)(\beta)$, where $Q$ is a map from
$\mathbb{R}^{n_1}$ and $\alpha,\beta$ are vectors in $\mathbb{R}^{n_1}$, to denote $\sum_{i,j}(D_iD_jQ)\alpha_i\beta_j$.

We now want to apply Feynman-Kac technique to deal with $\text{Ric}\nabla^HF$ as in Hsu~\cite{Hsu}. Let $M_t$ be the solution of the matrix-valued differential equation
$$\diff M_t+\frac{1}{2}M_t\text{Ric}_{U_t}\diff t=0, M_0=I.$$
Using It\^o's formula, it is straightforward to verify that
$$\diff M_t\nabla^HF(Z_t, T-t)+M_t\varPhi(Z_t, T-t)\diff t$$
is a martingale, where
$$\varPhi = (1)+(2)+(3)+(4)+(5).$$
Taking the expectation with respect to $\me_z$, we have
\begin{align*}
\nabla^H\E_zF(Z_T)=&\E_zM_t\nabla^HF(Z_t, T-t)\\
                     &+\E_z\int_0^tM_s\varPhi(Z_s, T-s)\diff s.
\end{align*}

Integrating over $t$ from $0$ to $T$ and using integration by parts on the second term on the right side, we have
\begin{align*}\tag{1}\label{eqn:1}
T\nabla^H\E_zF(Z_T)&=\E_z\int_0^TM_t\nabla^HF(Z_t, T-t)\diff t\\
                             &+\E_z\int_0^TN_t\varPhi(Z_t, T-t)\diff t.
\end{align*}
Here $N_t = (T-t) M_t$.  With this last identity, we have accomplished the task we have set ourselves at the beginning of this section, namely, we have passed the gradient operator through the expectation to act on the bundle map $f$, here represented by its scalarization $F$. As we pointed out above, the next step is to remove the expectation by an integration by parts argument in the path space.

\section{Integration by parts}

The first term on the right side of (1) contains the gradient on $F$.  From the heat equation $\frac{\partial F}{\partial t}=LF$ we have by It\^o's formula,
$$\diff F(Z_t, T-t)=\sum_i(H_i+V_{1,i})F\diff W_t^i.$$
For any $ v\in \mathbb{R}^{m}$ and $r\in\{1,2,\ldots,n_2\}$, we have
\begin{align*}
\E_z\int_0^T&\langle M_t\nabla^HF_r, v\rangle\diff
t\\
&=\E_z\int_0^T\langle \nabla^HF_r, M_t^\dagger
v\rangle\diff t\\
 &=\E_z\int_0^T\langle \nabla^HF_r, \diff
 W_t\rangle\int_0^T\langle M_t^\dagger v, \diff W_t\rangle\\
 &=\E_z\left(F_r(Z_T)-\E_zF_r(Z_T)-\int_0^T\sum_iV_{1,i}F_r\diff W_t^i\right)\int_0^T\langle v, M_t\diff W_t\rangle\\
 &=\E_z\left(F_r(Z_T)\int_0^T\langle v, M_t\diff
 W_t\rangle\right)-\E_z\int_0^T\langle\sum_iV_{1,i}F_re_i, M_t^\dagger
 v\rangle\diff t\\
 &=\E_z\left(F_r(Z_T)\int_0^T\langle v, M_t\diff
 W_t\rangle\right)-\E_z\int_0^T\langle M_t\sum_iV_{1,i}F_re_i,
 v\rangle\diff t\\
\end{align*}
This can be written more compactly as
\begin{align*}\tag{2}\label{eqn:2}
\E_z\int_0^T&M_t\nabla^HF\diff t\\
            &=\E_z\left(F(Z_T)\int_0^TM_t\diff W_t\right)-\E_z\int_0^TM_t\sum_jV_{1,j}Fe_j\diff t.
\end{align*}
The left side above is just the first term on the right side of (1). Note that we have removed the gradient operator from $F$. The remaining differentiation on $F$ are in the fibre directions.

The dealings with the remaining terms in (1) and the second term on the right side of (2)
$$\E_z\int_0^TN_t\varPhi(Z_t, T_t)-M_t\sum_j(DF)(V_{1,j})e_j\diff t$$
are admittedly technical, but in the process we will single out the terms involving some vertical derivatives on which we can perform a second commutation operation.

We will adopt the following notations: $\bar{N}_t^j$ stands for the horizontal lift of $N_te_j$ and $K_{t,i}^{j(l)}$ for the $2$nd and $3$rd coordinates of $[\bar{N}^j_t, H_i]$. Note that $N_t$ is symmetric, therefore in explicit components we have
$$N_te_j=\sum_kN_{t,j}^ke_k,\qquad\bar{N}^j_t=\sum_kN_{t,k}^jH_k,\quad\text{\rm and} \quad K_{t,i}^{j(l)}=\sum_k\Omega_{ki}^{(l)}N_{t,k}^j.$$We h
Recall that $\varPhi = (1) + (2) +(3)+(4)+(5)$. We have
$$N_t\varPhi(Z_t, T-t)-M_t\sum_j(DF)(V_{1,j})e_j=-A_1+A_2+A_3$$
where
$$A_1=\sum_j\sum_i\left(K_{t,i}^{j(2)}(H_iF)+\frac{1}{2}(H_iK_{t,i}^{j(2)})F+K_{t,i}^{j(2)}(DF)(V_{1,i})\right)e_j$$
$$A_2=\sum_j\sum_i\left(\left(DH_iF+(D^2F)(V_{1,i})\right)
\left(\bar{N}^j_tV_{1,i}+K_{t,i}^{j(1)}Y_t\right)\right)e_j$$
$$A_3=\sum_jDF\left(B_{t,j}+(\bar{N}_t^jV_0)-\sum_iV_{1,i}M_{t,i}^j\right)e_j$$
\begin{align*}
B_{t,j}&=\sum_i\frac{1}{2}\left((H_iK_{t,i}^{j(1)})Y_t+(D\bar{N}_t^jV_{1,i})(V_{1,i})\right.\\
       &\left.+(DV_{1,i})(\bar{N}_t^jV_{1,i})+(DV_{1,i})(K^{j(1)}_{t,i}Y_t)+H_i\bar{N}_t^jV_{1,i}\right).\\
\end{align*}
Here as will be in the sequel, we have omitted the ubiquitous $(Z_t, T-t)$ for the ease of notation.
The rule for regrouping the terms here may appear to be mysterious. Basically everything related to the curvature on $E_2$ is included in $A_1$, and the remaining terms follow the rule that the second order terms of $F$ go to $A_2$, and the first order terms go to $A_3$.


We consider the term $A_1$. From
$$\diff F=\sum_i(H_iF+(DF)(V_{1,i}))\diff W_t^i$$
we have
\begin{align*}
A_1\diff t&=\sum_j\left(\sum_iK^{j(2)}_{t,i}(H_iF+(DF)(V_{1,i}))\right)\diff te_j\\
          &+\frac{1}{2}\sum_j\left(\sum_i (H_iK^{j(2)}_{t,i})F\right)\diff te_j\\
          &=\sum_j\left(\sum_iK^{j(2)}_{t,i}\diff W_t^i\right)\cdot\left(\diff F\right)e_j\\
          &+\frac{1}{2}\sum_j\left(\sum_i (\diff K^{j(2)}_{t,i})\cdot(\diff W_t^i)\right)Fe_j\\
          &=\sum_j\left(\diff G_t^{j(2)}\cdot\diff F\right)e_j+\sum_j (\diff G_t^{j(2)})Fe_j+\diff S_t\\
          &=\sum_j\diff (G_t^{j(2)}F)e_j+\diff S'_t
\end{align*}
Here $G_t^{j(l)}$ are matrix valued process defined by by
$$G_t^{j(l)}=\int_0^t\sum_iK^{j(l)}_{t,i}\circ\diff W_t^i,$$
and $S_t$ and $S'_t$ are martingales. Integrating and taking the expectation we have
$$\E_z\int_0^TA_1\diff t=\E_z\sum_j\left(G_T^{j(2)}F(Z_T)\right)e_j.$$


\begin{remark}
As we mentioned before this paragraph of calculation, the curvature on $E_2$ gives rise to a separate term in the formula, and this term turns out
to only rely on the $0$th order term of the terminal value of $F$.
\end{remark}

\section{Second commutation and d\'enouement}

Before starting the calculations of $A_2, A_3$,  it is perhaps helpful to point out that if the vector bundle $E_2$ flat, then $A_1$ vanishes, but $A_2$ and $A_3$ remain.
We therefore anticipate terms involving the derivative $DF$.
Every term in $A_2$ and $A_3$ includes the $0$th or $1$st order term of
$DF$. Following Hsu~\cite{Hsu1}, we carry out a second commutation computation. From the heat equation for $F$ we have
\begin{align*}
\frac{\partial DF}{\partial t}&=LDF+[D,L]F\\
                                      &=LDF+A_4+A_5
\end{align*}
where
$$A_4=\sum_i\left((D^2F)(V_{1,i})+DH_iF\right)(DV_{1,i}),$$
$$A_5=DF\sum_i\left(\frac{1}{2}\left((D^2V_{1,i})(V_{1,i})+(DV_{1,i})(DV_{1,i})
+DH_iV_{1,i}\right)+DV_0\right).$$
Here $(D^2P)Q$ stands for the row vector $(\sum_iP_{ij}Q_i)$. As a consequence we have
$$\diff(DF)(Z_t,T-t)=\sum_i\left((D^2F)(V_{1,i})+DH_iF\right)\diff W_t^i-A_4\diff t-A_5\diff t$$
We now use a method similar to the one used in dealing with $A_1$ to absorb $A_4$ and $A_5$.
%
%
%

for this purpose we introduce the processes $\mr^{n_1}$-valued processes $Y_0^j=0$ called the derived processes in Norris~\cite{Norris}. They are determined by the stochastic differential equation
$$\diff Y_t^j=C_{2, i}^j\partial W_t^i+C_3^j\partial t+C_{4,i}^j\partial W_t^i+C_5^j\partial t,$$
where
$$C_{2,i}^j=\bar{N}^j_tV_{1,i}+K_{t,i}^{j(1)}Y_t,$$
\begin{align*}
C_3^j&=(\bar{N}_t^jV_0)-\sum_iV_{1,i}M_{t,i}^j\\
&+\sum_i\frac{1}{2}\left((H_iK_{t,i}^{j(1)})Y_t+(D\bar{N}_t^jV_{1,i})(V_{1,i})\right.\\
&\left.+(DV_{1,i})(\bar{N}_t^jV_{1,i}+K^{j(1)}_{t,i}Y_t)+H_i\bar{N}_t^jV_{1,i}\right),
\end{align*}
$$C_{4,i}^j=(DV_{1,i})Y_t^j$$
$$C_5^j=\sum_i\left(\frac{1}{2}\left((D^2V_{1,i})(V_{1,i})+(DV_{1,i})(DV_{1,i})
+DH_iV_{1,i}\right)+DV_0\right)Y_t^j.$$
The choices of these coefficients are made so that $A_4$ and $A_5$ can be eliminated explicitly from the final representation. They have the effect that
\begin{align*}\tag{3}\label{eqn:3}
&\qquad\left(N_t\varPhi(Z_t, T-t)-M_t\sum_j(DF)(V_{1,i})e_j\right)\diff t\\
&=-\sum_j\diff\left(G_t^{j(2)}F(Z_t)\right)e_j+\sum_j\diff\left(DF(Z_t)Y_t^j\right)e_j+\diff S_{1,t},
\end{align*}
where $S_{1,t}$ is a martingale.
In terms of Stratonovich integrals, the equation for the derived process can be written more compactly as
\begin{align*}
\diff Y_t^j&=-\sum_i(M_{t,i}^j)(V_{1,i})\diff t\\
     &\qquad+(\diff G_t^{j(1)})Y_t+(DV_0)Y_t^j\diff t+\bar{N}_t^jV_0\diff t\\
     &\qquad +\sum_i(DV_{1,i})(Y_t^j)\circ\diff W_t^i+\sum_i(\bar{N}_t^jV_{1,i})\circ\diff W_t^i
\end{align*}

Taking expectations in (3) and using (1) and (2) we finally obtain the desired probabilistic representation of the gradient for a solution of the heat equation in vector bundles.

\begin{theorem}  We have
$$T\nabla^H\E_z(F(Z_T))=\E_z\left[F(Z_T)\int_0^TM_t\diff W_t-G_TF(Z_T)+DF(Z_T)Y_T\right].$$
\end{theorem}

\begin{remark} In the setting of Norris~\cite{Norris}, his representation contains one extra term not present in our work. We have not worked through his work diligently to locate the source of this discrepancy.
\end{remark}

\end{document}